\documentclass[12pt]{article}
\usepackage{mathrsfs}
\usepackage{amsmath,amsfonts,amssymb,rotating,amsthm}
\usepackage{hyperref}
\usepackage{array}
\usepackage{breqn}
\usepackage{color}
\usepackage{fancyheadings}
\usepackage{mathptmx}
\usepackage{dsfont}
\usepackage{caption}
\usepackage{subcaption}
\usepackage{tikz}
\usetikzlibrary{trees}

\tikzstyle{every node}=[circle,inner sep=1pt,fill=white!60]

\tikzstyle{tn}=[shape=circle, draw, color=black!70]
\tikzstyle{marke}=[shape=circle,minimum size=0.2cm, draw,blue]

\allowdisplaybreaks
\def\qed{\nopagebreak\hfill{\rule{4pt}{7pt}}}
\def\proof{\noindent {\it{Proof.} \hskip 2pt}}
\newtheorem{thm}{Theorem}[section]
\newtheorem{deff}[thm]{Definition}

\newtheorem{prop}[thm]{Proposition}
\newtheorem{cor}[thm]{Corollary}

\theoremstyle{remark}

\numberwithin{equation}{section}

\usetikzlibrary{trees}

\tikzstyle{every node}=[circle,inner sep=1pt,fill=white!60]
\tikzstyle{tn}=[shape=circle, draw, color=black!70]
\tikzstyle{tn1}=[shape=circle, draw]
\tikzstyle{marke}=[shape=circle,minimum size=0.1cm, draw,blue]

\newcommand\gen{{\rm Gen}}

\begin{document}

\begin{center}
{\large\bf The Dumont Ansatz for the Eulerian Polynomials,

 Peak Polynomials and  Derivative Polynomials}

\vskip 6mm

\vskip 6mm

William Y.C. Chen$^1$ and Amy M. Fu$^2$

\vskip 3mm

$^{1}$Center for Applied Mathematics\\
Tianjin University\\
Tianjin 300072, P.R. China

\vskip 3mm

$^{2}$School of Mathematics\\
Shanghai University of Finance and Economics\\
Shanghai 200433, P.R. China

\vskip 3mm

Emails: {\tt $^1$chenyc@tju.edu.cn, $^{2}$fu.mei@mail.shufe.edu.cn}

\vskip 6mm

    {\bf Abstract}
\end{center}

We observe that three context-free grammars of Dumont can be brought to a
common ground, via the idea of
transformations of grammars, proposed by Ma-Ma-Yeh. Then we develop a
unified perspective to investigate several combinatorial objects in
connection with the bivariate  Eulerian polynomials. We call this
approach the Dumont ansatz.
 As applications, we provide grammatical
treatments, in the spirit of the symbolic method,
 of relations on the Springer numbers, the Euler numbers,
 the three kinds of peak polynomials, an identity of
 Petersen, and the two kinds of
derivative polynomials, introduced by Knuth-Buckholtz and  Carlitz-Scoville,
and later by Hoffman in a broader context. We obtain a convolution formula on the left peak polynomials,
leading to the  Gessel formula.
In this framework,  we come to the combinatorial interpretations of the derivative polynomials
due to Josuat-Verg\`es.

\noindent{\bf Keywords:} Dumont ansatz,  Eulerian polynomials, peak polynomials, derivative polynomials, context-free grammars

\noindent{\bf AMS Classification:} 05A15, 05A19

\section{Introduction}

The theme of this work is to present
a unified approach
to several enumeration
problems in connection with
the classical Eulerian polynomials
via a formal calculus based on context-free grammars.
We call this approach the Dumont ansatz because
it is largely built on the three grammars of
Dumont related to the Eulerian polynomials.
In some sense, we may say that the Dumont ansatz is a handful of
the grammars of  Dumont, reinforced by the idea of Ma-Ma-Yeh
\cite{Ma-Ma-Yeh-2019} concerning transformations of grammars.

The grammar of Dumont for the Eulerian polynomials reads
\begin{equation}\label{g-e}
G=\{ x \rightarrow xy, \;\; y \rightarrow xy\}.
\end{equation}
Let $D$ be the formal derivative with respect to the
above grammar $G$. For $n\geq 0$, the
bivariate Eulerian polynomial $A_n(x,y)$ is
defined by
\[ A_n(x, y) = D^n(y). \]

Recall that a grammar $G$ on a set  $V=\{x_1, x_2, \ldots\}$
of variables
is defined to be a set of substitution rules mapping each
variable $x_i$ to a Laurent polynomial $F_i(x_1, x_2, \ldots)$ on $V$, and the
formal derivative $D$ with respect to $G$ can be expressed as a differential
operator
\[ D=\sum_i F_i(x_1, x_2, \ldots)  {\partial \over \partial x_i} .\]
The generating function of a Laurent polynomial $f$ on $V$
is defined by
\begin{equation} \label{gen-f}
\gen(f, t) = \sum_{n=0}^\infty D^n(f) {t^n \over n!}.
\end{equation}
If $g$ is also a Laurent polynomial on $V$, then  $D$ satisfies the
product rule
\begin{equation}
     D(fg) =f D(g) + D(f) g .
\end{equation}
In general, $D$ obeys the Leibniz
rule,  i.e., for $n\geq 0$,
\begin{equation}\label{Leibniz-fg}
    D^n(fg) = \sum_{k=0}^n {n \choose k} D^k(f) D^{n-k}(g),
\end{equation}
or equivalently, the following multiplicative property holds,
\begin{equation} \label{gen-fg}
    \gen(fg, t) = \gen(f,t) \,\gen(g,t).
\end{equation}

For the above grammar $G$ in \eqref{g-e}, we have
\[ D= xy \left( {\partial \over \partial x} + {\partial \over \partial y} \right). \]
Setting $y=1$, the bivariate Eulerian polynomials $A_n(x,y)$
reduce to the Eulerian polynomials $A_n(x)$. The generating function
of $A_n(x)$ is given by
\begin{equation}
\label{Gen-A-n-x}
\sum_{n\geq 0} A_n(x) {t^n \over n!} = {1-x\over 1-x e^{(1-x)t}}.
\end{equation}
See, for example,  \cite{Petersen-2015}.
The above relation is equivalent to the bivariate version
\begin{equation} \label{g-a-x-y}
\sum_{n=0}^\infty A_n(x,y){t^n \over n!} =
{y-x \over 1 - xy^{-1} e^{(y-x)t} }.
\end{equation}

A grammatical derivation of (\ref{g-a-x-y})
is given in \cite{Chen-Fu-2022}.
The following form of the generating function of $A_n(x,y)$
for $n\geq 1$  is due to
Carlitz-Scoville \cite{Carlitz-Scoville-1972}, and it is
equivalent to the expression in the univariate case,
\begin{equation} \label{g-a-1}
\sum_{n=1}^\infty A_n(x,y) {t^n \over n!}
= xy \,{ e^{xt} - e^{yt}\over  xe^{yt} -y e^{xt}}.
\end{equation}

It is a well-known fact
due to Foata-Sch\"utzenberger  \cite{Foata-Schutzenberger-1970}
that the Eulerian polynomials $A_n(x)$ have the $\gamma$-expansion
with nonnegative coefficients,
\[
   A_n(x) = \sum_{k=1}^{\lfloor (n+1)/2 \rfloor}
    \gamma_{ n,k}  \,  x^k(1+x)^{n+1-2k},
\]
where $\gamma_{n,k}$ are nonnegative.

Let $u=xy$ and $v=x+y$. Then the grammar for the Eulerian polynomials
takes the form
\[ G=\{u \rightarrow uv, \;\; v\rightarrow 2u\}\]
and $A_n(x,y)$ can be expressed as $D^n(u)$.
For $1\leq n \leq 6$, the $\gamma$-expansions of $A_n(x,y)$ are as follows,
\begin{eqnarray*}
A_1(x,y) & = & u, \\[6pt]
A_2(x,y) & = & uv, \\[6pt]
A_3(x,y) & = &  u  v^2+2 u^2, \\[6pt]
A_4(x,y) & = & u v^3+8 u^2 v,\\[6pt]
A_5(x,y) & = &  u v^4+22 u^2 v^2+16 u^3, \\[6pt]
A_6(x,y) & = & u v^5+52  u^2 v^3+136 u^3 v.
\end{eqnarray*}

While the above grammar serves the purpose for the computation of
the $\gamma$-coefficients of the Eulerian polynomials, for the reason
that will be seen later, there is an advantage  to make the substitutions
\begin{equation} \label{u-xy-sub}
   u=xy, \;\; 2v = x+y.
\end{equation}
In this way, the transformed grammar becomes
\begin{equation} \label{g-d-2}
    G=\{u \rightarrow 2uv, \;\; v\rightarrow  u\},
\end{equation}
which is exactly the grammar given by Dumont for
increasing binary trees.
Then we define the bivariate Dumont polynomials, denoted $D_n(u,v)$,
in terms of the formal derivative of the grammar (\ref{g-d-2}),
or equivalently, in terms of increasing binary trees.

Under the substitutions in (\ref{u-xy-sub}), we can express $x$ and $y$ in terms
of $u$ and $v$, to wit,
\begin{equation}\label{xy-f-uv}
    x = v+ \sqrt{v^2 -u}  , \;\; y= v- \sqrt{v^2-u}.
\end{equation}

It is remarkable that Dumont also discovered an
analogous grammar for 0-1-2 increasing trees, that is,
\[ G=\{u \rightarrow uv, \;\; v\rightarrow  u\}.\]
Despite the striking resemblance of the aforementioned three grammars,
they have been playing their own roles without supporting
each other. Thanks to the idea of transformations of grammars, due to Ma-Ma-Yeh
\cite{Ma-Ma-Yeh-2019},
we recognize that these three grammars can be brought to a
common ground in the name of the Dumont ansatz.

Applying the strategies of the Dumont ansatz to the  peak polynomials
 and the derivative polynomials,
we illustrate how to make connections to
the Eulerian polynomials. A left peak of a permutation is
also called an exterior peak or a peak. There are other two
relevant statistics, the number  of interior peaks  and the
number of left-right peaks (or outer peaks). We give
grammatical labelings, called the $M$-labeling and the
$W$-labeling, so that we can bring the three kinds of
peak polynomials to the test ground of the Dumont ansatz.

It should be stressed that once we have a grammar
on file, it is often not hard to find a combinatorial
structure of a recursive nature such as increasing plane
trees or increasing binary trees, as an interpretation
of the corresponding polynomials. This means that
a grammar can be instrumental in search for appropriate
combinatorial structures.

The derivative polynomials $P_n(x)$ and $Q_n(x)$ for the tangent
and the secant were introduced by Knuth-Buckholtz \cite{Knuth-Buckholtz-1967}
in their studies of the tangent, Euler and Bernoulli numbers.
They were studied later by Carlitz-Scoville \cite{Carlitz-Scoville-1972} and
Hoffman \cite{Hoffman-1995, Hoffman-1999} in broader contexts, see also
\cite{Ma-2012}.
When evaluated
at $x=1$, the derivative polynomials $Q_n(x)$ turn out to be the
Springer numbers.
Using the Dumont ansatz, we quickly get the combinatorial interpretations
of the derivative polynomials established  by Josuat-Verg\`es
\cite{Josuat-2014}.  Moreover, we see how the
the derivative polynomials $P_n(x)$ are related to
the Eulerian polynomial $A_n(x)$.

The Dumont ansatz not only provides a mechanism to unify
many known results, but also offers a rigorous platform to exploit
the grammatical calculus, in the spirit of the symbolic
method, in the course of proving and discovering
combinatorial identities.
For example, we found it possible to give a derivation of an identity of
Petersen by using the grammatical calculus. We also
obtain a convolution identity on the
left peak polynomials, which can be used to derive the
 formula of Gessel on the generating function of the left
 peak polynomials
\cite[Sequence A008971]{OEIS}, see also  \cite{Zhuang-2016}.

\section{The Dumont ansatz }

In this section, we first give an overview of three grammars
of Dumont, for the Eulerian polynomials, increasing binary trees, and
the Andr\'e polynomials. We see that
these three grammars share the same nature and the
corresponding generating functions can be deduced from each other
via a change of variables. We use an intermediate structure
as a unified model to deal with the relations among the generating
functions in the family, and we call this approach the Dumont
ansatz. Roughly speaking, the idea behind the Dumont ansatz is to
establish connections among combinatorial polynomials by means of
transformations of grammars.

\subsection{The grammar for the Eulerian polynomials}

As described in the Introduction,
the bivariate Eulerian polynomials $A_n(x,y)$ can be
 expressed as $D^n(y)$, where $D$ is the formal derivative
 with respect to the grammar
 \begin{equation}\label{g-a}
 G=\{ x \rightarrow xy, \;\;
 y \rightarrow xy\}.
 \end{equation}
The  bivariate polynomials  $A_n(x,y)$
can be expressed in terms of the numbers of descents
and ascents of permutations, or in terms of complete increasing binary
trees, which are called the Gessel trees in \cite{Chen-Fu-Yan-2022}.

For $n=0$, we define $A_0(x,y)=y$, the first few values of $A_n(x,y)$ are given below,
\begin{eqnarray*}
A_1(x,y) & = & xy, \\[6pt]
A_2(x,y) & = & xy^2+x^2y, \\[6pt]
A_3(x,y) & = & xy^3+4x^2y^2+x^3y, \\[6pt]
A_4(x,y) & = & xy^4 + 11 x^2y^3 + 11 x^3 y^2 + x^4y,\\[6pt]
A_5(x,y) & = & xy^5+ 26x^2y^4 +66x^3y^3 +26 x^4y^2 + x^5y, \\[6pt]
A_6(x,y) & = & xy^6+57x^2y^5+302x^3y^4+ 302 x^4y^3 + 57 x^5y^2 + x^6y.
\end{eqnarray*}
For more information about context-free grammars for
combinatorial polynomials and Eulerian polynomials,
see \cite{Chen-Fu-2017, Chen-Fu-2022, Dumont-1996}.
The journey of the Dumont ansatz starts with the above grammar $G$ and the Eulerian polynomials $A_n(x,y)$.

\subsection{The grammar for increasing binary trees}

The second grammar of Dumont we will be concerned with is
\begin{equation} \label{g-uv}
    G=\{ u \rightarrow 2uv, \;\; v \rightarrow u\}.
\end{equation}
Here we deliberately choose to use the variables $u$ and $v$
rather than $x$ and $y$ as in \cite{Dumont-1996}, because
it is related to the grammar (\ref{g-a}) for the Eulerian polynomials
via the following substitutions
\begin{equation}\label{sub-uv-xy}
    u=xy, \;\; 2v=x+y.
\end{equation}
Since
\[ D(u)=D(xy)=xy(x+y) = 2uv\]
and
\[  D(v) = D(x+y)/2 = xy = u,\]
the grammar $G$ in (\ref{g-a}) is transformed into the
grammar $G$ in (\ref{g-uv}). This transformation implies that
the Eulerian polynomials $A_n(x)$ are $\gamma$-positive, as
observed by Ma-Ma-Yeh \cite{Ma-Ma-Yeh-2019}.

\begin{figure}[!ht]
\begin{center}
\begin{tikzpicture}[scale=0.9]
\node [tn,label=180:$1$]{}[grow=down]
	%[sibling distance=16mm,level distance=10mm]
    child[grow=231] {node [tn,label=180:{$4(v)$}](four){}
     %[sibling distance=14mm,level distance=13mm]
            child {node [tn,label=180:{$6$}](two){}
        %[sibling distance=14mm,level distance=13mm]
        child[grow=230] {node [tn,label=-90:{$9(u)$}](three){}}
        child[grow=310] {node [tn,label=-90:{$7(u)$}](three){}}
        }}
    child[grow=310] {node [tn,label=0:{$2$}](two){}
        %[sibling distance=14mm,level distance=10mm]
        child[grow=231] {node [tn,label=0:{$3(v)$}](three){}
        %[sibling distance=14mm,level distance=13mm]
    child [grow=305]{node [tn,label=-90:{$5(u)$}](){}
    }}
        child[grow=310] {node [tn,label=0:{$8(u)$}](three){}}
        };

\end{tikzpicture}
\end{center}
\caption{An increasing binary tree with the $(u,v)$-labeling.}
\label{fibt-labeled}
\end{figure}
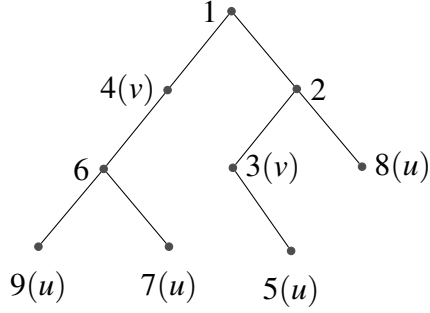

Dumont \cite{Dumont-1996} showed that the
grammar $G$ in (\ref{g-uv}) gives a weighted counting
of increasing binary trees $T$ on $[n]$ if we label
a leaf by $u$ and a degree one vertex by $v$, where the degree of
a vertex in a binary tree is referred to the number of its children.
Then we define the weight of $T$ as the product of
the labels associated with $T$. The labels of the
increasing binary tree in Figure \ref{fibt-labeled}  are shown in
parentheses. This way of labeling an increasing binary tree is called
the $(u,v)$-labeling.

The initiative of the Dumont ansatz grew out of the realization that
the underlying combinatorial structures of the original grammar
and the transformed grammar are essentially the same. The substitutions
of variables are reflected by different labeling schemes.
It should be mentioned that the idea of using a change
of variables to compute the $\gamma$-coefficients of the
Eulerian polynomials has
appeared in the work of Chow \cite{Chow-2008}.

\begin{deff}
For $n\geq 0$, the Dumont polynomial $D_n(u,v)$ is referred to
the polynomial $D^n(v)$, where $D$ is the formal derivative with
respect to the grammar $G$ in (\ref{g-uv}).
\end{deff}

The following theorem is due to Dumont \cite{Dumont-1996}.

\begin{thm}[Dumont] \label{thm-D-uv}
For $n\geq 1$, $D_n(u,v)$ equals the sum of weights
of all increasing binary trees on $[n]$ endowed with the $(u,v)$-labeling.
\end{thm}

The first few values of
$D_n(u,v)$ are listed  below,
\begin{eqnarray*}
D_0(u,v) & = & v, \\[6pt]
D_1(u,v) & = & u, \\[6pt]
D_2(u,v) & = & 2uv, \\[6pt]
D_3(u,v) & = &  4 u v^2+2 u^2, \\[6pt]
D_4(u,v) & = &  8 u v^3+16 u^2 v,\\[6pt]
D_5(u,v) & = &  16 u v^4+88 u^2 v^2+16 u^3, \\[6pt]
D_6(u,v) & = & 32 u v^5+416 u^2 v^3+272 u^3 v.
\end{eqnarray*}

\subsection{0-1-2 increasing plane trees}

In the study of the $\gamma$-coefficients of the
Eulerian polynomials, 0-1-2 increasing plane trees
arise as an underlying structure for
the combinatorial interpretation.
In accordance with the grammar $G$ in (\ref{g-uv}),
we need an alternative
labeling scheme.

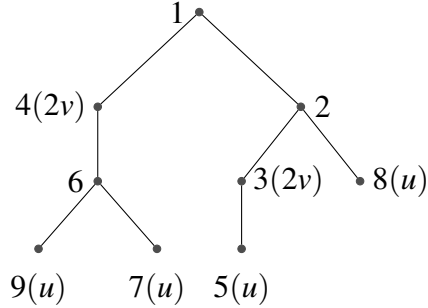
\begin{figure}[!ht]
\begin{center}
\begin{tikzpicture}[scale=0.9]
\node [tn,label=180:$1$]{}[grow=down]
	[sibling distance=30mm,level distance=14mm]
    child {node [tn,label=180:{$4(2v)$}](four){}
    	[sibling distance=14mm,level distance=11mm]
            child {node [tn,label=180:{$6$}](two){}
        [sibling distance=17.5mm,level distance=10mm]
        child {node [tn,label=-90:{$9(u)$}](three){}}
        child {node [tn,label=-90:{$7(u)$}](three){}}
        }}
    child {node [tn,label=0:{$2$}](two){}
        [sibling distance=17.5mm,level distance=11mm]
        child {node [tn,label=0:{$3(2v)$}](three){}
        [sibling distance=17.5mm,level distance=10mm]
        child {node [tn,label=-90:{$5(u)$}](nine){}}}
        child {node [tn,label=0:{$8(u)$}](three){}}
        };
\end{tikzpicture}
\end{center}
\caption{A 0-1-2 increasing plane tree with the $(u,2v)$-labeling.}
\label{fip-labeled}
\end{figure}

Let $n\geq 1$, and let $T$ be a 0-1-2
increasing plane tree on $[n]$. A $(u,2v)$-labeling of $T$ is
referred to labeling every leaf by $u$ and labeling every
degree one   vertex by $2v$.
Given a $(u,2v)$-labeling, the weight of $T$ is
defined to be the product of labels associated with $T$.
Figure \ref{fip-labeled} depicts a 0-1-2 increasing plane tree with
the $(u, 2v)$-labeling.

Using the $(u,2v)$-labeling of 0-1-2 increasing plane trees,
the theorem of Dumont concerning the grammar (\ref{g-uv}) can be
reformulated as follows.
%  which is implicitly known to OEIS \cite[A]{OEIS}.

\begin{thm}[Dumont]
For $n\geq 1$,
the Dumont polynomial $D_n(u,v)$ equals the sum of weights
of all 0-1-2 increasing plane trees on $[n]$ with the $(u,2v)$-labeling.
\end{thm}

It is evident that the $(u, 2v)$-labeling for 0-1-2 increasing plane trees
is equivalent to the $(u,v)$-labeling for increasing binary trees.
Thus,
for $n\geq 1$,
\begin{equation}\label{d-n-a-n}
D_n(u,v)= A_n(x,y),
\end{equation}
where $(u,v)$ and $(x,y)$ are related by (\ref{sub-uv-xy}), i.e.,
$u=xy$ and $2v=x+y$.

\subsection{The Andr\'e polynomials}

Dumont \cite{Dumont-1996} showed that the following grammar
$$
G= \{u \rightarrow uv, \;\;v\rightarrow u\}
$$
can be used to generate the bivariate Andr\'e polynomials $E_n(x,y)$,
which are defined in terms of 0-1-2 increasing trees as follows,
 $$
E_n(u,v)=\sum_{T}u^{f_0(T)}v^{f_1(T)},
$$
where the sum ranges over all 0-1-2 increasing trees $T$ on $[n]$, and
$f_i(T)$ denotes the numbers of degree $i$ vertices of $T$ for $i=0, 1, 2$.
For $n\geq 1$ and
$x=y=1$, $E_n(x,y)$
reduces to the Euler number $E_n$, that is, the number of alternating permutations
on $[n]$.

Define $E_0(u,v)=1$.
The first few values of $E_n(u,v)$ are given below,
\begin{eqnarray*}
E_1(u,v) & = & u, \\[6pt]
E_2(u,v) & = & u v, \\[6pt]
E_3(u,v) & = &   u v^2+u^2, \\[6pt]
E_4(u,v) & = & u v^3+4 u^2 v ,\\[6pt]
E_5(u,v) & = & u v^4+11 u^2 v^2+4 u^3 , \\[6pt]
E_6(u,v) & = & u v^5+26 u^2 v^3+34 u^3 v.
\end{eqnarray*}

Notice that the bivariate version of the Andr\'e polynomials
can be recovered from the one variable version. Let $E_n(u) = E_n(u, 1)$.
Then for $n\geq 1$,
\begin{equation}\label{E-n-u-v-2}
       E_n(u, v) = v^{n+1} E_n\left( {u \over v^2} \right).
\end{equation}
The generating function of the Andr\'e polynomials $E_n(x)$ was
obtained by Foata-Sch\"utzenberger \cite{FS-1973}. An alternative proof can be found in  Foata-Han \cite{FH-2001}. A derivation  utilizing
the grammar of Dumont was given in \cite{Chen-Fu-2017}.
The following  is manifest since
both are generated by essentially the same grammar.

\begin{thm} \label{thm-DE}
For $n\geq 0$, we have
\begin{equation}
D_n(2u,v) = 2^n E_n(u,v).
\end{equation}
\end{thm}

\proof
Clearly, ordering the two children of a degree two vertex in a
0-1-2 increasing tree is equivalent to assigning the number two
to this vertex as a label, so that
\begin{equation}\label{d-n-u-v-T}
D_n(u,v) =  \sum_{L} 2^{f_2(T)} u^{f_{0}(T)} (2v)^{f_1(T)},
\end{equation}
where the sum ranges over 0-1-2 increasing trees on $[n]$.
It follows that
\[ D_n(2u, v) = \sum_T 2^{f_0(T) + f_1(T) + f_2(T)}
                      u^{f_0(T)} v^{f_1(T)} = 2^n \sum_T
                      u^{f_0(T)} v^{f_1(T)},
 \]
where $T$ has the same range as in (\ref{d-n-u-v-T}),
whereupon the theorem is proved.  \qed

Now we see that  the Andr\'e polynomials can be expressed in terms of
the Eulerian polynomials, see \cite[A094503]{OEIS}.
For $n\geq 0$, we have
\begin{equation} \label{enan}
   2^n E_n(u,v) = A_n(x, y),
\end{equation}
where $x$ and $y$ are determined by $xy=2u$ and $x+y=2v$.

For $u=1$ and $v=1$, (\ref{enan}) becomes the known
identity on the Euler numbers,
\begin{equation} \label{enani}
    E_n = {A_n(i) \over (1+i)^{n-1}},
\end{equation}
where $n\geq 1$ and $i = \sqrt{-1}$, see \cite[Sequence A000111]{OEIS}.

\section{The peak polynomials}

The objective of this section is to demonstrate that
the peak polynomials, all of the three kinds,
fall into the framework of the Dumont ansatz.
First of all, let us now get the notation straight.
Our proposal is to employ symbols that are
meaningful and yet easy to remember. It turns out that
the $LMW$-notation seems to be a sensible choice.
As for left peaks (exterior peaks), the letter $L$ looks like
having a peak on the left, and so we use $L(n,k)$ to denote
the number of permutations of $[n]$ with $k$ left peaks.
Accordingly, we use $L_n(x)$ and $L_n(x,y)$ to denote
the one variable and bivariate  left peak polynomials, respectively, that is,
\begin{equation} \label{tnkx-2}
    L_n(x) = \sum_{k=0}^{\lfloor n /2 \rfloor} L(n,k) x^{k},
\end{equation}
and
\begin{equation} \label{tnkx}
    L_n(x,y) = \sum_{k=0}^{\lfloor n /2 \rfloor} L(n,k) x^{2k+1}
        y^{n-2k}.
\end{equation}
The first few values of $L_n(x,y)$ are given
below,
\begin{eqnarray*}
L_0(x,y) &=& x, \\[6pt]
L_1(x,y) & = &  xy , \\[6pt]
L_2(x,y) & = &  xy^2+x^3, \\[6pt]
L_3(x,y) & = &  xy^3 +5 x^3y , \\[6pt]
L_4(x,y) & = &  xy^4+18 x^3y^2 +5 x^5,\\[6pt]
L_5(x,y) & = &  xy^5 +58x^3y^3 +61 x^5y  , \\[6pt]
L_6(x,y) & = &  xy^6 +179x^3y^4+479x^5y^2+61 x^7 .
\end{eqnarray*}

To be more specific, let $n \geq 1$ and
let $\sigma=\sigma_1\sigma_2\cdots \sigma_n$  be a permutation of $[n]$.
We
assume that $\sigma_0=\sigma_{n+1}=0$. Then an index $i$  is said
to be a left peak if $ 1 \leq i <n$ and $\sigma_{i-1} < \sigma_i > \sigma_{i+1}$,
or an interior peak  if $ 1 < i <n$ and $\sigma_{i-1} < \sigma_i > \sigma_{i+1}$,
or a left-right peak  if $ 1 \leq i \leq n$ and $\sigma_{i-1} < \sigma_i > \sigma_{i+1}$,

Next, we choose the letter $M$ for the case of interior peaks,
because the two peaks in $M$ bear a striking resemblance to interior peaks.  Therefore,
we shall use $M(n,k)$ to denote the number of permutations of $[n]$ with
$k$ interior peaks.
For $n\geq 1$, the interior peak
polynomials are defined by
\begin{equation}
    M_n(x) = \sum_{k=0}^{\lfloor (n-1)/2 \rfloor } M(n,k) x^{k}  .
\end{equation}

In the case of left-right peaks or outer peaks,
the letter $W$
signifies three left-right peaks including
the two at both ends, and so we use
$W(n,k)$ to denote the number of permutations of $[n]$
with $k$ left-right peaks. We move on to define
\begin{equation}
    W_n(x) = \sum_{k=1}^{\lfloor (n+1) /2 \rfloor} W(n,k) x^{k} .
\end{equation}

Note that various notations for the peak polynomials and their
coefficients
 have appeared in the literature, see, for example,
\cite{Carlitz-1974, Chow-Ma-2014, Kitaev-2007, Ma-2012, MFMY-2022, Petersen-2007},
while they do not necessarily mean the same as in here.
The bivariate versions of $M_n(x)$ and $W_n(x)$ are crucial
as far as the grammars are concerned, which are defined by
\begin{eqnarray} \label{def-mn-xy}
 M_n(x,y) & = &  \sum_{k=0} ^{\lfloor (n-1)/2\rfloor} M(n,k) x^{2k+2} y ^{n-2k-1}, \\[9pt]
 W_n(x,y) & = &  \sum_{k=1} ^{\lfloor (n+1)/2\rfloor} W(n,k)x^{2k} y ^{n-2k+1}.
\end{eqnarray}
In fact, there is a reason to express $W_n(x,y)$ as
\begin{equation} \label{wnxy-2}
W_n(x,y)   =    \sum_{k=0} ^{\lfloor (n-1)/2\rfloor} W(n,k+1)x^{2k+2} y ^{n-2k-1}.
\end{equation}

The first few values of $W_n(x)$ are given below,
\begin{eqnarray*}
W_0(x,y)  & = & y ,   \\[6pt]
W_1(x,y)  & = & x^2 ,   \\[6pt]
W_2(x,y)  & = & 2 x^2 y ,   \\[6pt]
W_3(x,y)  & = & 4  x^2y^2+2x^4 ,   \\[6pt]
W_4(x,y)  & = & 8 x^2y^3+16x^4  y ,   \\[6pt]
W_5(x,y)  & = & 16 x^2y^4+88x^4y^2+16 x^6 ,   \\[6pt]
W_6(x,y)  & = & 32 x^2y^5+416x^4y^3+272x^6y.
\end{eqnarray*}

It can be seen from the above table that
the polynomials $W_n(x)$ have the same coefficients as
the Dumont polynomials $D_n(u,v)$.

\subsection{Connection to the Gessel formula}

Adopting the Dumont ansatz,
we obtain a convolution formula connecting the left peak polynomials with the
Dumont polynomials, which yields  the Gessel formula on the
generating function of $L_n(x)$,  i.e.,
\[ L(x,t) = \sum_{n=0}^\infty L_n(x ) { t^n \over n!}. \]
Similarly, let
\[ M(x,t) = \sum_{n=0}^\infty M_n(x ) { t^n \over n!}. \]
Note that by way of
recurrence relations,  David-Barton \cite{David-Barton-1962}
established partial differential equations
on $L(x,t)$ and  $M(x,t)$ and found solutions in differential forms.

\begin{thm}[David-Barton]
We have
\begin{eqnarray}
 2x(1-x) {\partial L(x,t) \over \partial x}  + (xt-1)
 {\partial L(x,t) \over \partial t} + L(x,t) +1 & = & 0, \\[9pt]
 2x(1-x) {\partial M(x,t) \over \partial x}  + (xt-1)
 {\partial M(x,t) \over \partial t} + M(x,t) + \sqrt{x} & = & 0.
\end{eqnarray}
Put
\[ z = t \sqrt{1-x} + \log \left( { \sqrt{x} \over 1+ \sqrt{1-x}} \right). \]
Then we have
\begin{eqnarray} \label{ptxt}
  {\partial L(x,t) \over \partial x} & = & {1 \over 2}
      \left(  {1 \over \cosh(z) -1} + { 1 \over \cosh (z) +1} \right), \\[9pt]
  {\partial M(x,t) \over \partial t} & = &{1 \over 2}
      \left(  {1 \over \cosh(z) -1} - { 1 \over \cosh (z) +1} \right). \label{pmxt}
\end{eqnarray}
\end{thm}

Gessel   \cite[Sequence A008971]{OEIS} obtained the
following explicit  formula for $L(x,t)$.

\begin{thm}[Gessel] We have
\begin{equation}  \label{Gessel-F}
L(x,t)
     =\frac{\sqrt{1-x}}{\sqrt{1-x}\cosh{(t\sqrt{1-x})}-\sinh{t\sqrt{1-x}}}.
\end{equation}
\end{thm}

As brought up by Stanley \cite{Stanley-V-1}, an explicit expression for
$M(x,t)$ can be deduced from Equation (\ref{pmxt}) of David-Barton. An extension
to a more general enumeration problem was given by Carlitz-Scoville
\cite{Carlitz-Scoville-1972}.

Our point of departure is the
following grammar
\begin{equation} \label{g-peak}
G=\{  x\rightarrow xy, \quad y \rightarrow x^2 \},
\end{equation}
independently found by Chen-Fu \cite{Chen-Fu-2017} and Ma \cite{Ma-2012}.
Let $D$ be the formal derivative of the grammar $G$.
It has been shown that the left peak polynomials $L_n(x,y)$ can be generated by
the grammar $G$, i.e.,  for $n\geq 0$,
\begin{equation}
    L_n(x,y) = D^n(x).
\end{equation}

Setting $u=x^2$ and $v=y$, the grammar in (\ref{g-peak}) takes the form
\begin{equation} \label{g-peak-2}
G=\{  u \rightarrow 2uv, \;\; v \rightarrow u\},
\end{equation}
which turns out to be exactly a grammar of the Dumont ansatz.
This transformation enables us to establish a convolution identity
on $L_n(x,y)$.

Let $D$ be the formal derivative
with respect to the grammars in (\ref{g-peak}) and (\ref{g-peak-2}). Bear in mind
that there is no ambiguity
because the substitution rules act on distinct variables.
Let us consider
the polynomials $D^{n}(v)$.
There are two ways to look at $D^{n}(v)$.
On one hand,   $D_{n}(u, v)$ can be considered as
a polynomial in $u,v$, which equals the Dumont polynomial $D_{n}(u,v)$.
On the other hand, $D^n(v)$ can be treated as a polynomial in $x,y$.

Since
\[ D^{n+1}(y) = D^{n }(x^2),\]
we obtain  the following convolution identity.

\begin{thm} \label{thm-anc}
For $n\geq 0$,
\begin{equation}
    D_{n+1}(x^2,y) =  x^2 \sum_{k=0} ^{n }
    { n   \choose k} L_{k}(x,y)  L_{n -k}(x,y). \label{anc}
\end{equation}
\end{thm}

The above formula gives rise to a relation on
the generating function of $L_n(x,y)$. Let $L(x,y,t)$ be the generating function of  $L_n(x,y)$.

\begin{cor} Let $A(\bar{x},\bar{y},t)$ be the generating function of
the Eulerian polynomials $A_n(\bar{x},\bar{y})$, and let $A'(\bar{x},\bar{y},t)$
denote the differentiation with respect to $t$. Then we have
\begin{equation} \label{apuvt}
 A'(\bar{x},\bar{y},t)= L^2(x,y,t),
\end{equation}
where $\bar{x}=y+\sqrt{y^2-x^2}$ and $\bar{y}=y-\sqrt{y^2-x^2}$.
\end{cor}

Let us illustrate how to compute  $L(x,y,t)$ from $A(\bar{x},\bar{y},t)$.
Invoking (\ref{g-a-1}), that is,
\[ A(\bar{x},\bar{y},t) = \bar{x}\bar{y} \,{ e^{\bar{x}t} - e^{\bar{y}t}\over  \bar{x}e^{\bar{y}t} -\bar{y} e^{\bar{x}t}}.
\]
we find that
\begin{equation}\label{dauv}
    A'(\bar{x},\bar{y},t) =\bar{x}\bar{y} e^{ (\bar{x}+\bar{y})t}  { (\bar{x}-\bar{y})^2 \over (\bar{x}  e^{\bar{y} t}-\bar{y} e^{\bar{x} t})^2}.
\end{equation}
In order to connect $A'(\bar{x},\bar{y},t)$ with  $L(x,y,t)$,
it is necessary to express the Dumont polynomial $D_n(u,v)$ in terms of
the Eulerian polynomial $A_n(x,y)$. Substituting $u=x^2$ and $v=y$ into
(\ref{xy-f-uv}), we find that
\begin{equation}
    A_n(y+\sqrt{y^2-x^2},  y-\sqrt{y^2-x^2}) = D_n(x^2, y).
\end{equation}
Plugging
\[ \bar{x}=y+\sqrt{y^2-x^2}, \quad \bar{y}=y-\sqrt{y^2-x^2}\]
into (\ref{dauv}), a routine calculation shows that
\[
   A'(y+\sqrt{y^2-x^2}, y-\sqrt{y^2-x^2},t) =
   \left(\frac{ xy \sqrt{y^2-x^2}}
   {\sqrt{y^2-x^2}\cosh{(t\sqrt{y^2-x^2})}-y\sinh{t\sqrt{y^2-x^2}}} \right)^2,
\]
which,  together with (\ref{apuvt}), gives
\begin{equation} \label{txyt}
    L(x,y,t)= \frac{xy\sqrt{y^2-x^2}}{\sqrt{y^2-x^2}\cosh{(t\sqrt{y^2-x})}-y\sinh{t\sqrt{y^2-x^2}}}.
\end{equation}
Noting that $\sqrt{x}L_n(x) =  L_n(\sqrt{x},1)$,
so that $L(x,t) =  \sqrt{x}L(\sqrt{x}, 1 ,t)$.
By (\ref{txyt}), we obtain (\ref{Gessel-F}).

 Petersen \cite{Petersen-Thesis}
established a
relation between $L_n(x)$ and the Eulerian polynomials,
see \cite{Ma-2012}, to wit,
\begin{equation} \label{PI}
    L_n \left( {4x \over (1+x)^2 }\right)
    = {1 \over (1+x)^n}
       \sum_{k=0}^n {n \choose k} (1-x)^{n-k} 2^k
           A_k(x)   .
\end{equation}

To seek a grammatical understanding of the above identity,
we find a rather simple transformation
forging a bridge between $L_n(x,y)$ and $A_n(x,y)$.

\begin{thm}\label{thm35} Setting $x=\sqrt{\bar{x}\bar{y}}$ and
$y=(\bar{x} + \bar{y})/2$,
the grammar \[ G=\{\bar{x}\rightarrow \bar{x}\bar{y}, \;\; \bar{y} \rightarrow \bar{x}\bar{y}\}\]
for the Eulerian polynomials is transformed into the grammar
\[
G=\{x \rightarrow xy, \quad y\rightarrow x^2\}
\]
for the left peak polynomials.
\end{thm}

\proof By definition, we have
$$
D(x)=D(\sqrt{\bar{x}\bar{y}})= {D(\bar{x}\bar{y}) \over 2 \sqrt{\bar{x}\bar{y}}}
= \sqrt{\bar{x}\bar{y}} \,\frac{\bar{x}+\bar{y}}{2} = xy
$$
and
$$
D(y)=D\left({\bar{x}+\bar{y} \over 2}\right) = \bar{x}\bar{y}=y^2,
$$
as requested. \qed

It should be mentioned that while we usually
consider Laurent polynomials for the action of a
formal derivative. But as noted in \cite{Chen-1993},
we are not confined to Laurent polynomials. Indeed, taking
 square root is not an issue at all. To help understand the formal
 derivative of $\sqrt{\bar{x}\bar{y}}$, we may set $z =\sqrt{\bar{x}\bar{y}} $ and then
 apply the product rule
 \[ D(z^2)= D(\bar{x} \bar{y}) = \bar{x} \bar{y} (\bar{x}+\bar{y}) = 2 z D(z) \]
 to obtain $D(z)$.

 \noindent{\it Grammatical Proof of Petersen's Identity (\ref{PI}).}
Recall that
$$
D^n(x)=L_n(x,y)=\sum_{k=0}^{[n/2]} L(n,k)x^{2k+1} y^{n-2k}.
$$
In view of the grammar transformation given in
Theorem \ref{thm35}, we see that
$$
D^n(\sqrt{\bar{x}\bar{y}})=L_n\left(\sqrt{\bar{x}\bar{y}},\frac{\bar{x}+\bar{y}}{2}\right).
$$
On the other hand, by the Leibniz rule, we have
\begin{eqnarray*}
D^n(\sqrt{\bar{x}\bar{y}})&=&
D^n\left(\frac{y}{\sqrt{\bar{x}^{-1}\bar{y}}}\right) \\[9pt]
& = & \sum_{k=0}^n{n \choose k}D^k(\bar{x})D^{n-k}(\sqrt{\bar{x}\bar{y}^{-1}})\\[9pt]
&=&\sqrt{\bar{x}\bar{y}^{-1}}\sum_{k=0}^n {n \choose k} A_k(\bar{x},\bar{y})\frac{(\bar{y}-\bar{x})^{n-k}}{2^{n-k}},
\end{eqnarray*}
where we have made use of the fact that for $k\geq 0$,
\begin{equation}
   D^k(\bar{x}\bar{y}^{-1})=\bar{x}\bar{y}^{-1}(\bar{y}-\bar{x})^k,
 \end{equation}
as observed in \cite{Chen-Fu-2017}. It follows that
\begin{equation}\label{TA-1}
L_n\left(\sqrt{\bar{x}\bar{y}},\frac{\bar{x}+\bar{y}}{2}\right)
=\sqrt{\bar{x}\bar{y}^{-1}}\sum_{k=0}^n {n \choose k} A_k(\bar{x},\bar{y})\frac{(\bar{y}-\bar{x})^{n-k}}{2^{n-k}}.
\end{equation}
Setting $\bar{y}=1$ and replacing $\bar{x}$ by $x$ yields
\begin{eqnarray}\label{TA-2}
L_n\left(\sqrt{x},\frac{1+x}{2}\right)  &=&
\sqrt{x}\sum_{k=0}^n {n \choose k} A_k(x)\frac{(1-x)^{n-k}}{2^{n-k}}\nonumber\\[9pt]
&=&\frac{\sqrt{x}}{2^n}\sum_{k=0}^n {n \choose k} 2^kA_k(x)(1-x)^{n-k}.
\end{eqnarray}
On the other hand,
\begin{eqnarray}\label{TA-3}
L_n\left(\sqrt{x},\frac{1+x}{2}\right) &=&
\sum_{k=0}^{\lfloor n /2 \rfloor} L(n,k) (\sqrt{x})^{2k+1}\frac{(1+x)^{n-2k}}{2^{n-2k}}\nonumber\\[9pt]
&=&\frac{\sqrt{x}(1+x)^n}{2^n}\sum_{k=0}^{\lfloor n /2 \rfloor} L(n,k) \frac{(4x)^k}{(1+x)^{2k}}.\end{eqnarray}
Comparing (\ref{TA-2}) and (\ref{TA-3}),  we arrive at \eqref{PI}.
\qed

\subsection{The grammatical labelings}

A grammatical labeling of permutations was given in \cite{Chen-Fu-2017}
to produce the left peak polynomials. Similar tactics can be applied
to the other two kinds of peak polynomials.

The labeling for $L_n(x,y)$, called the $L$-labeling,
can be described as follows. Let $\sigma$ be a
permutation of $[n]$. We patch a zero to $\sigma$ at both ends so that there are
$n+1$ positions between two adjacent elements, and these are the possible positions
to insert $n+1$ in $\sigma$ to generate a permutation of $[n+1]$.
The $L$-labeling of $\sigma$ is meant to label the last position by $x$, label
the two positions next to any left peak by $x$, and label the remaining positions by $y$.
Notice that the $L$-labeling is an equivalent representation of the labeling given
in \cite{Chen-Fu-2017}.
Below is an example,
\[ 3 \;1 \;4 \;5 \;6 \;2 \;   \stackrel{L}{\longrightarrow} \;
  0 \; x \; 3 \; x  \;1 \; y \;4 \; y \; 5 \;x \;6 \; x\; 2 \; x \; 0 .\]

The procedure of generating the above permutation, along with the
$L$-labelings and the corresponding substitution rules,
is displayed in the table below,
\begin{center}
\begin{tabular}{|c|l|l|l|}
\hline
$\;\;n\;\;$  &  The $L$-labeling & Weight  &
Rule $\quad\;\;$ \\ \hline
1 & $ 0\;  y \;1 \;x \; 0 $ &  $  xy $ & $ x \rightarrow xy $\\ \hline
2 & $ 0 \; y \;1 \; y\;2 \; x \; 0  $ &  $  xy^2  $ & $ x \rightarrow x y  $\\ \hline
3 & $ 0 \;x \; 3 \; x \;1 \; y \;2 \;  x\;0 $ &  $  x^3y $ & $ y \rightarrow x^2 $\\ \hline
4 & $ 0 \;x \; 3 \; x \;1 \; x\; \; 4 \; x \;2 \;  x\;0 $ &  $  x^5 $ & $ y \rightarrow x^2$\\ \hline
5 & $ 0 \;x \; 3 \; x \;1 \; y \; 4 \; x \;5 \; x \;2 \;  x\;0 $ &  $  x^5y $
& $ x \rightarrow xy$
\\ \hline
6 & $ 0 \;x \; 3 \; x \;1 \; y \; 4 \; y\; 5 \; x \; 6 \; x \;2 \;  x\;0 \qquad $
&  $  x^5y^2 $
& $ x \rightarrow xy$\\ \hline
\end{tabular}
\end{center}
\vskip 6pt

The $M$-labeling of $\sigma$ is supposed
to label the two positions at both ends by $x$,
label the two positions next to any interior peak by $x$, and label the remaining
positions by $y$. Below is an example,
\[ 3 \;1 \;4 \;5\;6 \;2  \; \stackrel{M}{\longrightarrow} \;
  0 \; x \; 3 \; y  \;1 \; y \;4 \; y \; 5 \;x \;6 \;x \;2 \; x \; 0 .\]
The procedure of generating the above permutation along with the
$M$-labelings is displayed in the table below,
\begin{center}
\begin{tabular}{|c|l|l|l|}
\hline
$\;\;n\;\;$  &  The $M$-labeling & Weight  &
Rule $\quad\;\;$ \\ \hline
1 & $ 0\;  x \;1 \;x \; 0 $ &  $  x^2 $ & $ y \rightarrow x^2 $\\ \hline
2 & $ 0 \; x \;1 \; y\;2 \; x \; 0  $ &  $  x^2 y  $ & $ x \rightarrow x y  $\\ \hline
3 & $ 0 \;x \; 3 \; y \;1 \; y \;2 \;  x\;0 $ &  $  x^2y^2 $ & $ x \rightarrow xy $\\ \hline
4 & $ 0 \;x \; 3 \; y \;1 \; x\; \; 4 \; x \;2 \;  x\;0 $ &  $  x^4y $ & $ y \rightarrow x^2$\\ \hline
5 & $ 0 \;x \; 3 \; y \;1 \; y \; 4 \; x \;5 \; x \;2 \;  x\;0 $ &  $  x^4y^2 $
& $ x \rightarrow xy$
\\ \hline
6 & $ 0 \;x \; 3 \; y \;1 \; y \; 4 \;y \; 5 \; x \; 6 \; x \;2 \;  x\;0 \qquad $
&  $  x^4y^3 $
& $ x \rightarrow xy$\\ \hline
\end{tabular}
\end{center}
\vskip 6pt

The following relation was derived by Ma \cite{Ma-2012} via recurrence
relations.

\begin{thm}
For $n\geq 1$, we have
\begin{equation} \label{dnm}
    D^n(y) = M_n(x,y).
\end{equation}
\end{thm}

Indeed, it is only a matter of formality to check that the $M$-labeling
justifies the above conclusion.
Let us now turn to the $W$-labeling.  For a permutation $\sigma$ of
$[n]$, a zero is patched at both ends. Consider the positions between
two adjacent elements of $\sigma$, where $\sigma_0=\sigma_{n+1}=0$. For
any element $1\leq i \leq n$, if $\sigma_i$ is a left-right peak, then
label the positions next to $i$ by $x$, and label the rest of the positions by
$y$.

For example, the $W$-labeling of the permutation $1$ is $0\;x\;1\;x\;0$, and below
is a permutation accompanied by its $W$-labeling,
\[ 3 \;1 \;4 \;5 \; 6\;2 \; \stackrel{W}{\longrightarrow} \;
 0 \; x \; 3 \; x  \;1 \; y \;4 \; y \; 5\; x\; 6 \;x \;2 \; y \; 0 .\]

By examining the change of labels when inserting the element $n+1$ to
a permutation on $[n]$, we are led to the following interpretation of the
grammatical expansion of $D^n(y)$.
The same reasoning as for the grammatical generation of
the polynomials $L_n(x,y)$ given in \cite{Chen-Fu-2017} applies to $W_n(x,y)$.
Thus the justification with full rigor will not be repeated here.
Instead, we shall give an example. The process of generating the preceding
permutation
is described in the following table,
\begin{center}
\begin{tabular}{|c|l|l|l|}
\hline
$\;\;n\;\;$  &  The $W$-labeling & Weight  &
Rule $\quad\;\;$ \\ \hline
1 & $ 0\;  x \;1 \;x \; 0 $ &  $  x^2 $ & $ y \rightarrow x^2 $\\ \hline
2 & $ 0\; y  \;1 \; x \;2 \;  x \; 0  $ &  $  x^2 y  $ & $ x \rightarrow x y  $\\ \hline
3 & $ 0   \; x \; 3 \; x   \;1 \; x  \;2 \;x  \;0 $ &  $  x^4  $ & $ y \rightarrow x^2 $\\ \hline
4 & $ 0 \; x\; 3 \;   x \;1 \;  x \; 4 \; x  \;2 \;  y\;0 $ &  $  x^4y $ & $ x \rightarrow xy$\\ \hline
5 & $ 0 \; x \; 3 \; x  \;1 \; y  \; 4 \;  x \;5 \; x  \;2 \; y  \;0 $ &  $  x^4y^2 $
& $ x \rightarrow xy$
\\ \hline
6 & $ 0 \; x \; 3 \;  x \;1 \;  y \; 4 \; y \; 5 \;  x \; 6 \; x  \;2 \; y  \;0 \qquad $
&  $  x^4y^3 $
& $ x \rightarrow xy$\\
\hline
\end{tabular}
\end{center}
\vskip 6pt

\begin{thm}
For $n\geq 1$, we have
\begin{equation} \label{dnw}
    D^n(y) = W_n(x,y).
\end{equation}
\end{thm}

Now that $W_n(x) = xM_n(x)$,
by comparing (\ref{def-mn-xy}) with (\ref{wnxy-2}),
 we retrieve the following combinatorial property, see Zhuang \cite{Zhuang-2016}.

\begin{thm}
For $n\geq 1$ and $0\leq k \leq \lfloor (n-1)/2 \rfloor$, the number
of permutations of $[n]$ with $k+1$ left-right peaks
equals the number of permutations of $[n]$ with $k$
interior peaks, that is,
\begin{equation}
    W(n,k+1) = M(n,k).
\end{equation}
\end{thm}

Applying the Dumont ansatz, we get
the following relation, where $D_{n}(u,v)$ are
the Dumont polynomials.

\begin{thm} \label{thm-dm}
For $n\geq 0$,
\begin{equation}
    D_{n}(u,v) = M_n(x,y),
\end{equation}
where $u=x^2$ and $v=y$.
\end{thm}

Using the theory of enriched $P$-partitions,
Stembridge \cite{Stembridge-1997} showed that for $n\geq 1$,
\begin{equation}
   x M_n\left( {4 x \over (1+x)^2} \right) = {2^{n-1} \over (1+x)^{n-1}}
    \,A_n(x),
\end{equation}
see also, \cite{Ma-2012, Petersen-2015}.
In fact,  Theorem
\ref{thm-dm} spells out how the these polynomials (in the bivariate forms)
are related via a change of variables.

\subsection{Convolution formulas}

From the point of view of grammars,
one can convert the Leibniz formulas into
convolution formulas for combinatorial polynomials,
and we shall exemplify how this is the case for peak polynomials.

Returning to $M_n(x)$ and $W_n(x)$, we see that for $n\geq 1$,
\begin{equation}\label{MW-1}
x M_n(x) = W_n(x),
\end{equation}
since $ x^2M_n(x^2) =  M_n(x,1)$ and $ W_n(x^2) =  W_n(x,1)$.

In consideration of the initial values, we have no choice but to impose that
\begin{equation}
    M_0(x) =x^{-1} \quad \mbox{and} \quad W_0(x)= 1.
\end{equation}

Since $D(x)= xy$, we have for $n\geq 0$,
\[ D^{n+1}(x) = D^n(xy) .\]
and so
  the following convolution formula is immediate.

\begin{thm}
For $n\geq 0$, we have
\begin{equation} \label{CTM-1}
    L_{n+1}(x, y) = \sum_{k=0}^n {n \choose k} L_{k}(x,y) M_{n-k}(x,y),
\end{equation}
or equivalently,
\begin{equation} \label{CTM-2}
    L_{n+1}(x) = x \sum_{k=0}^n {n \choose k} L_{k}(x ) M_{n-k}(x).
\end{equation}
\end{thm}

Since $D(x^2)=D(y^2) = 2x^2y$, we deduce that for $n\geq 1$,
\[ D^{n} (x^2) =D^n(y^2) ,  \]
which can be translated into the following convolution identity.

\begin{thm} \label{thm-CTM}
For $n\geq 1$, we have
\begin{equation} \label{CTM}
   \sum_{k=0}^n {n \choose k} L_k(x,y)L_{n-k}(x,y)
   = \sum_{k=0}^n {n \choose k} M_{k}(x,y) M_{n-k}(x,y).
   \end{equation}
\end{thm}

Given that $D(x^2)=D(y^2)$, we find that
for $n\geq 1$,
\[D^{n+1}(y) = D^{n}(x^2) = D^n(y^2).\]
Therefore, we reach a convolution formula for $M_n(x,y)$, which
is reminiscent of that for the
derivative polynomials $P_n(x)$. This is no surprise
because the polynomials $M_n(x,y)$ and $P_n(x)$ admit
the same grammatical structure.

\begin{thm}\label{thm-CMXY}  For $n\geq 1$, we have
\begin{equation} \label{CMXY}
   M_{n+1}(x,y)
   = \sum_{k=0}^n {n \choose k} M_{k}(x,y) M_{n-k}(x,y).
   \end{equation}
\end{thm}

With the understanding that $M_0(x)=x^{-1}$ and $M_1(x)=1$,
for $n\geq 1$, (\ref{CMXY}) may be replaced by
\begin{equation} \label{CMXY2}
   M_{n+1}(x)
   = x \sum_{k=0}^n {n \choose k} M_{k}(x) M_{n-k}(x).
   \end{equation}

In the same vein, for $n\geq 1$,
the convolution formula (\ref{CTM}) can be recast as
\begin{equation} \label{CTMa}
    \sum_{k=0}^n {n \choose k} L_k(x)L_{n-k}(x)
   = \sum_{k=0}^n {n \choose k} M_{k}(x) M_{n-k}(x).
   \end{equation}

 Ma \cite{Ma-2012} considered the polynomials
$R_n(x,y)= D^n(x+y)$, where $n\geq 0$. Since
$$
D^{n+1}(x+y) = D^n (x(x+y)),
$$
we get the following convolution formula at once.

\begin{thm}
For $n\geq 0$, we have
\begin{equation}
    R_{n+1}(x,y)= \sum_{k=0}^n {n \choose k} L_{k}(x,y) R_{n-k}(x,y).
\end{equation}
\end{thm}

    The convolution of $L_n(x)$
appeared in Ma-Yeh \cite{Ma-Yeh-2016}. However, the above convolution identity
for $L_n(x)$ went unnoticed,
even though there is no barrier to deduce it in the context.

\section{The derivative polynomials}

The derivative polynomials $P_n(x)$ and $Q_n(x)$ for
the tangent and the secant were  introduced by Knuth-Buckholtz
\cite{Knuth-Buckholtz-1967}
for the computation of tangent, Euler and  Bernoulli numbers.
They were later studied by
Carlitz-Scoville
\cite{Carlitz-Scoville-1972}.
 Hoffman \cite{Hoffman-1995} defined the
 derivative polynomials in a more general context.
  Recall that
\[ (\tan(x))' =   \tan^2(x) + 1.\]
Define $P_n(x)$ by
\begin{equation} \label{d-p}
{{\rm d} ^n \over {\rm d} x^n} \tan(x) = P_n(\tan(x)).
 \end{equation}

 The derivative polynomials $P_n(x)$ are listed as the Sequence A008293 in OEIS, and
 $P_n(0)$ are the tangent numbers.
  As pointed out
 by Dumont \cite{Dumont-1995}, the numbers $P_n(1)$ go back to Euler,
 and they are listed as Sequence A000831 in OEIS \cite{OEIS} with
 initial values
 \[
1, 2, 4, 16, 80, 512, 3904, 34816, 354560, 4063232, \ldots .\]

 In the same vein,  the derivative
 polynomials $Q_n(x)$  for the secant are defined by
 \begin{equation} \label{s-p}
{{\rm d} ^n \over {\rm d} x^n} \sec(x) = Q_n(\tan(x)) \sec(x).
 \end{equation}

For $x=0$, $Q_n(x)$ are the secent numbers, that is,
$Q_n(0)=0$ for $n$ odd, and $Q_n(0)=E_n$ for $n$ even.
The numbers $Q_n(1)$
 are  called the Springer numbers or the generalized
  Euler numbers, denoted by
  $S_n$.

The derivative polynomials satisfy the recurrence relations for $n\geq 0$,
\begin{eqnarray}
    P_{n+1}(x) & = & (1+x^2) {{\rm d} \over {\rm d} x} P_n(x),
      \label{p-n-rec} \\[9pt]
    Q_{n+1}(x) & = &
 (1+x^2) {{\rm d} \over {\rm d} x} Q_{n }(x) + x Q_{n }(x).
\label{q-n-rec}
\end{eqnarray}

Note that the generating functions of $P_n(x)$ and $Q_n(x)$ were given by
Hoffman \cite{Hoffman-1995}.

\begin{thm}[Hoffman] We have
\begin{eqnarray}
\sum_{n=0}^\infty  P_n(x) {t^n \over n!} & = &
{x + \tan (t) \over 1-x \tan(t)}, \label{p-n-z-g}
\\[9pt]
\sum_{n=0}^\infty  Q_n(x) {t^n \over n!} & = &  {1 \over \cos(t) - x \sin(t)}. \label{q-n-z-g}
\end{eqnarray}
\end{thm}

Notice that the generating function of $Q_n(x)$ tells that
$Q_n(1)$ coincides with the Springer number $S_n$ for any $n$.

\subsection{The grammar}

Whereas the following theorem is merely a
paraphrase of the recursive formulas for $P_n(x)$ and $Q_n(x)$,
it lends a different perspective to look at the derivative polynomials,
as well as a channel to the Dumont ansatz.

\begin{thm} Let
\begin{equation} \label{g-q-x}
G =\{ a \rightarrow a x, \;\; x\rightarrow 1+ x^2 \},
\end{equation}
and let $D$ be the formal derivative with respect to $G$.
 Then for $n\geq 0$,
 \begin{eqnarray}
 D^n(x) & = & P_n(x), \\[6pt]
 D^n(a) & = & a Q_n(x).
 \end{eqnarray}
\end{thm}

For $n=0$, we have $P_0(x)=x$ and $Q_0(x)=1$. The first
few values of $P_n(x)$  and $Q_n(x)$ are given below,
 \begin{eqnarray*}
 P_1(x) & = & 1+x^2 , \\[6pt]
 P_2(x) & = & 2 x + 2 x^3, \\[6pt]
 P_3(x) & = & 2+8 x^2 + 6 x^4, \\[6pt]
 P_4(x) & = & 16 x +40 x^3 +  24 x^5, \\[6pt]
 P_5(x) & = & 16+136 x^2+240 x^4+ 120 x^6 ,\\[6pt]
 P_6(x) & = & 272 x +1232 x^3 +1680 x^5 + 720 x^7,
 \end{eqnarray*}
and
\begin{eqnarray*}
Q_1(x) & = & x,\\[6pt]
Q_2(x) & = & 1+x^2,\\[6pt]
Q_3(x) & = &  5x+6x^3,\\[6pt]
Q_4(x) & = & 5+28 x^2+24 x^4,\\[6pt]
Q_5(x) & = & 61 x+180 x^3+120 x^5, \\[6pt]
Q_6(x) & = & 61+662 x^2+1320 x^4+720 x^6.
\end{eqnarray*}

\subsection{The generating functions}

 We present a proof of the generating functions of $P_n(x)$ and
 $Q_n(x)$ to
demonstrate the rigor and the efficiency of the grammatical calculus.
In the formal setting, the generating functions of $P_n(x)$ and
$Q_n(x)$ can be expressed as follows.

\begin{thm} We have
\begin{eqnarray}
{\rm Gen} (a,t) & = & {a \over \cos (t) -x \sin(t)}, \label{gut}\\[9pt]
{\rm Gen}(x,t)& = &  {  x \cos (t) +  \sin (t)  \over  \cos (t) -  x \sin (t)}  .\label{gxt}
\end{eqnarray}
\end{thm}

\noindent{\it Grammatical Proof.} Since
\begin{equation}
{\rm Gen}(a, t) =  {1 \over {\rm Gen}(a^{-1}, t)} \label{g-u-t}
\end{equation}
and
\begin{equation}
{\rm Gen}(x,t) =  {{\rm Gen}(a^{-1}x,t) \over {\rm Gen}(a^{-1},t)}, \label{g-x-t}
\end{equation}
we proceed to compute  $D^n(a^{-1})$ and
$D^n(a ^{-1}x)$. Clearly,
\[ D(a^{-1}) = - a^{-2} D(a)= - a^{-2}(ax) =-a^{-1} x,\]
so we get
\[ D^2(a^{-1}) = D(-a^{-1}x) = a^{-1}x^2 - a^{-1}(1+x^2) = -a^{-1},\]
and hence
\[ D^3(a^{-1}) = a^{-1} x.\]
It follows that for $n\geq 0$,
\begin{eqnarray}
D^{2n}(a^{-1}) & = & (-1)^n a^{-1}, \label{d-2n}\\[6pt]
D^{2n+1}(a^{-1}) & = & (-1)^{n+1}a^{-1}x. \label{d-2n+1}
\end{eqnarray}
Consequently,
\begin{eqnarray*}
{\rm Gen}(a^{-1},t) &  = & \sum_{n=0}^\infty D^n(a^{-1}){t^n \over n!} \\[9pt]
 & = &   \sum_{n=0}^\infty (-1)^n a^{-1} {t^{2n} \over (2n)!} +
  \sum_{n=0}^\infty (-1)^{n+1} a^{-1}x {t^{2n+1} \over (2n+1)!} \\[9pt]
  & = & a^{-1} \cos (t) - a^{-1} x \sin (t).
\end{eqnarray*}
This proves (\ref{gut}).

In order to compute ${\rm Gen}(a^{-1}x,t)$, observe that (\ref{d-2n}) and (\ref{d-2n+1}) can
be alternatively expressed as
\begin{eqnarray}
D^{2n}(a^{-1}x) & = & (-1)^n a^{-1}x, \label{d-2n-a}\\[6pt]
D^{2n+1}(a^{-1}x) & = & (-1)^{n } a^{-1}, \label{d-2n+1-b}
\end{eqnarray}
from which we infer that
\begin{equation}\label{g-u-1-x}
    {\rm Gen}(a^{-1}x,t) = a^{-1} x \cos (t) + a^{-1} \sin(t).
\end{equation}
According to (\ref{g-x-t}), we reach (\ref{gxt}).
This completes the proof. \qed

\subsection{Convolution formulas}

To show how the grammatical calculus can be performed, we consider the
following convolution identities of Hoffman \cite{Hoffman-1995}. For $n\geq 1$,
\begin{eqnarray}
  P_{n+1}(x) & = & \sum_{k=0}^n {n\choose k} P_k(x) P_{n-k}(x), \label{PC} \\[9pt]
  Q_{n+1} (x) & = & \sum_{k=0}^n {n\choose k} P_k(x) Q_{n-k}(x). \label{QC}
\end{eqnarray}

Since $D^{n+1}(x) = D^n(1+x^2)=D^n(x^2)$,
(\ref{PC}) immediately follows from the Leibniz rule.
To verify (\ref{QC}), we only need the condition $n\geq 0$. Note that
$D^{n+1}(a) = D^n(ax)$. Thus  (\ref{QC}) is a consequence of the
Leibniz rule.

The following convolution formula found by Ma-Fang-Mansour-Yeh  \cite{MFMY-2022} also  follows from the same line of reasoning. For $n\geq 0$,
\begin{equation}\label{MFMY-convolution}
P_{n+2}(x) = 2 \sum_{k=0}^n {n \choose k} P_k(x) P_{n+1-k}(x).
\end{equation}
Noting that
\[ D^2(x) = D(1+x^2) = D(x^2) = 2xD(x),\]
we have
\[ D^{n+2}(x) = 2 D^n(xD(x)).\]
Applying the Leibniz rule yields (\ref{MFMY-convolution}).

The following  convolution identity is  due to Hoffman \cite{Hoffman-1995},
which can be used to evaluate $Q_n(1)$.

\begin{thm}[Hoffman]
For $n\geq 0$,
\begin{equation} \label{Convolution-3}
    P_{n+1}(x) = (1+x^2) \sum_{k=0}^n {n \choose k} Q_k(x) Q_{n-k}(x).
\end{equation}
\end{thm}

There is  a grammatical derivation of  (\ref{Convolution-3}). We need the following
property.

\begin{prop}
We have
  \begin{equation} \label{dxu}
 D\left( {1+x^2 \over a^2}\right) =0.
 \end{equation}\end{prop}

 \proof
  Since
  \[ D(a^{-2}) =-2a^{-3}D(a) = -2a^{-2}x\]
  and
  \[D(1+x^2)=2xD(x)=2x(1+x^2),\]
  we see that
  \[ D\left( {1+x^2 \over a^2}  \right)
  = {-2x(1+x^2) \over a^2} + {2x(1+x^2)
  \over a^2} =0, \]
  as claimed. \qed

\noindent{\it Grammatical Proof of (\ref{Convolution-3}).}
It suffices to show that
\begin{equation} \label{Convolution-4}
    a^2 D^{n+1}(x) = (1+x^2) \sum_{k=0}^n {n \choose k} D^k(a)  D^{n-k}(a).
\end{equation}
In view of the Leibniz rule, the above relation takes the form
\begin{equation} \label{Convolution-5}
    a^2 D^{n+1}(x) = (1+x^2) D^n(a^2).
\end{equation}
But (\ref{dxu}) says that ${1+x^2 \over a^2}$ is
a constant with respect to $D$,  this implies that
\[ {1+x^2 \over a^2} D^n(a^2) = D^n(1+x^2) = D^{n+1}(x),\]
as required. \qed

\subsection{The Josuat-Verg\`es trees}

As will be seen, the grammar $G$
is a reflection of the recursive construction of a
structure introduced by Josuat-Verg\`es which serves as
the basis for combinatorial
interpretations of the derivative polynomials.

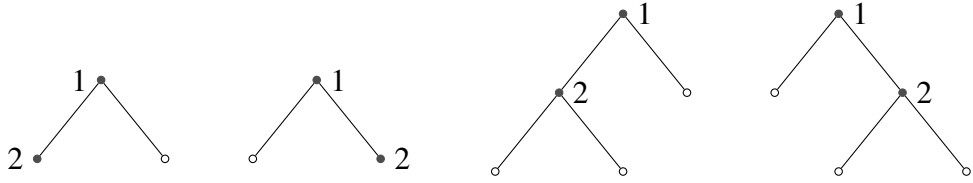
\begin{figure}[!ht]
\begin{center}
\begin{tikzpicture}[scale=0.9]
\node [tn,label=180:$1$]{}[grow=down]
	%[sibling distance=16mm,level distance=10mm]
    child[grow=231] {node [tn,label=180:{$2$}](one){}}
     child[grow=309] {node [tn1,label=0:{}](){}};
\end{tikzpicture}
\qquad
\begin{tikzpicture}[scale=0.9]
\node [tn,label=0:$1$]{}[grow=down]
	%[sibling distance=16mm,level distance=10mm]
    child[grow=231] {node [tn1,label=0:{}](one){}}
     child[grow=309] {node [tn,label=0:{$2$}](){}};
\end{tikzpicture}
\qquad
\begin{tikzpicture}[scale=0.9]
\node [tn,label=0:$1$]{}[grow=down]
	%[sibling distance=16mm,level distance=10mm]
    child[grow=231] {node [tn,label=0:{$2$}](one){}
     child[grow=231] {node [tn1,label=0:{}](one1){}}
      child[grow=309] {node [tn1,label=0:{}](one2){}}
    }
         child[grow=309] {node [tn1,label=0:{}](one1){}};
\end{tikzpicture}
\qquad
\begin{tikzpicture}[scale=0.9]
\node [tn,label=0:$1$]{}[grow=down]
	%[sibling distance=16mm,level distance=10mm]
    child[grow=231] {node [tn1,label=0:{}](one1){}}
    child[grow=309] {node [tn,label=0:{$2$}](one){}
     child[grow=231] {node [tn1,label=0:{}](one1){}}
      child[grow=309] {node [tn1,label=0:{}](one2){}}
    };
\end{tikzpicture}
\end{center}
\caption{The four Josuat-Verg\`es trees on $\{1, 2\}$.}
\label{jvt}
\end{figure}

Let us recall the definition.  Let $n\geq 1$.
A Josuat-Verg\`es tree is a complete increasing binary tree
possibly with  unlabeled leaves, called empty leaves.
More precisely, the vertices
from the root to any leaf are in the increasing order
regardless of the empty leaf, if any. For $n=0$, the
only Josuat-Verg\`es tree is the empty leaf.
The labeling stipulation of a Josuat-Verg\`es tree is simple:
each empty leaf is labeled by $x$.
For $n=2$, the four
Josuat-Verg\`es trees are listed in Figure \ref{jvt}, with the total weight
amounting to $P_2(x)=2x+2x^3$.

Now, it can be seen that the substitution rule $x\rightarrow 1+x^2$
corresponds to adding the element
$n+1$ to a Josuat-Verg\`es tree on $[n]$
by turning an empty leaf into a new vertex with no children
or an internal vertex  with
two empty children. This operation would result in a Josuat-Verg\`es tree on
$[n+1]$, and it offers a recipe to understand the following statement.

\begin{thm}
For $n\geq 1$, $P_n(x)$ equals the sum
of weights of all Josuat-Verg\`es trees on
$[n]$.
\end{thm}

\subsection{Exponentiation of the Dumont grammar}

Examining the partition argument for the Fa\`a de Bruno formula
in regard with  higher order derivatives of the composition of two functions,
or equivalently, the combinatorial interpretation for successive applications
of
the formal derivative of the grammar
\[ G= \{ f_i \rightarrow f_{i+1} g_1, \;\; g_j \rightarrow g_{j+1} \;|\, i = 0, 1, 2, \ldots, \; j=1,2, \ldots \}, \]
see \cite{Chen-1993},
we are led to a combinatorial interpretation of $Q_n(x)$
based on that of $P_n(x)$. The
proof is solely a repetition of the aforementioned testimony, and hence is omitted.

To be more specific, we consider the grammar
\begin{equation} \label{g-e-d}
    G=\{ a \rightarrow a    v, \;\;v \rightarrow u, \;\; u \rightarrow 2uv\}.
\end{equation}
Let $D$ denote the formal derivative with respect to $G$. The partition argument
shows that $D^n(a)$ corresponds to the exponential structure built on
the structure of the original grammar of Dumont. This makes it possible to
lift the Dumont ansatz to the exponential level, involving
forests of planted increasing binary trees.

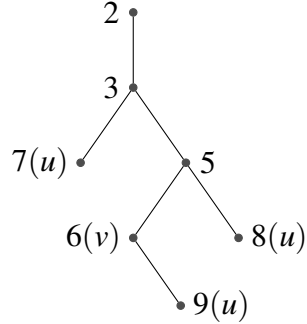
\begin{figure}[!ht]
\begin{center}
\begin{tikzpicture}
\node [tn,label=180:$2$]{}[grow=down]
	[sibling distance=14mm,level distance=10mm]
    child {node [tn,label=180:{$3$}](four){}
        child {node [tn,label=180:{$7(u)$}](one){}}
        child {node [tn,label=0:{$5$}](two){}
        [sibling distance=14mm,level distance=10mm]
        child {node [tn,label=180:{$6(v)$}](three){}
        [sibling distance=14mm,level distance=11mm]
    child [grow=305]{node [tn,label=0:{$9(u)$}](){}
    }}
        child {node [tn,label=0:{$8(u)$}](three){}}
        }
        };
\end{tikzpicture}
\end{center}
\caption{A planted increasing binary tree.}
\label{fpibt}
\end{figure}

Let us take up Theorem \ref{thm-D-uv} as an example. It is obvious that for
$n\geq 1$,
$D^n(v)$ equals the sum of weights of planted increasing binary trees
on $[n]$. A forest of planted
increasing
binary trees on $[n]$ can be naturally defined as a forest on $[n]$ in
the usual sense, with each component being a planted increasing binary tree.
The labeling schemes for planted binary trees can be carried over to forests.
However, we should pay special attention to a tree with a single vertex.
If a leaf   differs from the root, we label it with $u$. If a degree one vertex
is not the root, we label it with $v$. If the root is the only vertex
in the tree, then we label it with $v$ as it is the starting
label, otherwise the root is not labeled,
bearing in mind that there is at most one child of the root.
For example, Figure \ref{fpibt} depicts a planted increasing binary tree
with labels in parentheses.

The exponential counterpart of Theorem \ref{thm-D-uv}
can be described as follows.

\begin{thm}\label{thm-E-D-uv}
For $n\geq 1$, $D^n(a)$ equals the sum of weights
of all forests of planted increasing binary trees on $[n]$ with the $(u,v)$-labeling.
\end{thm}

In the framework of the Dumont ansatz, we are
ready to make a connection to
the  Josuat-Verg\`es forests on $[n]$, which are referred to as
plane rooted forests in \cite{Josuat-2014}. First, a planted
Josuat-Verg\`es tree on a nonempty subset $S$ of $[n]$ is
defined as an increasing rooted tree possibly with empty leaves
such that the root has only one child, which is allowed to be an
empty leaf, and the subtree of the root, if not empty, is a Josuat-Verg\`es
tree on the rest of the elements in $S$. For the special case  when $S$ contains
only one element, it is necessary to
note  that the unique planted Josuat-Verg\`es tree on $S$
consists of the root along with an empty leaf.

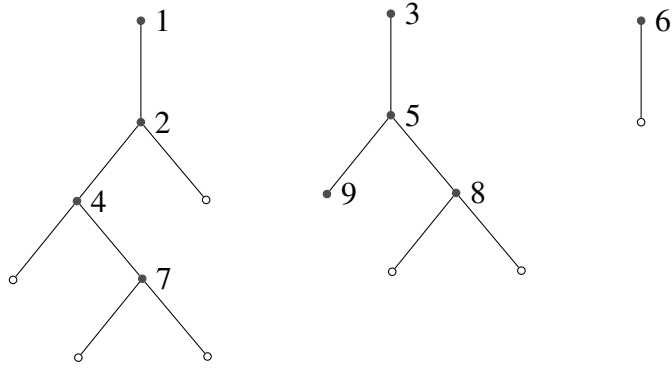
\begin{figure}[!ht]
\begin{center}
\begin{tikzpicture}[scale=0.9]
\node [tn,label=0:$1$]{}[grow=down]
	%[sibling distance=16mm,level distance=10mm]
    child {node [tn,label=0:{$2$}](two){}
    child[grow=231] {node [tn,label=0:{$4$}](four){}
     %[sibling distance=14mm,level distance=13mm]
     child[grow=231] {node [tn1,label=0:{}](five){}}
     child[grow=310] {node [tn,label=0:{$7$}](seven){}
     child[grow=231] {node [tn1,label=0:{}](){}}
     child[grow=310] {node [tn1,label=0:{}](){}}
     }}
     child[grow=310] {node [tn1,label=0:{}](){}}
     };
\end{tikzpicture}
\qquad\quad
\raisebox{6ex}{\begin{tikzpicture}[scale=0.9]
\node [tn,label=0:$3$]{}[grow=down]
	%[sibling distance=16mm,level distance=10mm]
    child {node [tn,label=0:{$5$}](two){}
    child[grow=231] {node [tn,label=0:{$9$}](four){}}
    child[grow=310] {node [tn,label=0:{$8$}](eight){}
    child[grow=231] {node [tn1,label=0:{}](){}}
    child[grow=310] {node [tn1,label=0:{}](){}}
    }
     %[sibling distance=14mm,level distance=13mm]
     };
\end{tikzpicture}}
\qquad\quad
\raisebox{16.5ex}
{\begin{tikzpicture}[scale=0.9]
\node [tn,label=0:$6$]{}[grow=down]
	%[sibling distance=16mm,level distance=10mm]
    child {node [tn1,label=0:{}](two){} };
\end{tikzpicture}}
\end{center}
\caption{A Josuat-Verg\`es forest.}
\label{jvf}
\end{figure}

Then a Josuat-Verg\`es
forest is defined as a forest on $[n]$ consisting of planted Josuat-Verg\`es trees.
The weight of a Josuat-Verg\`es forest is defined  to be the
product of the weights of all trees in the forest.  For example, Figure \ref{jvf} demonstrates a Josuat-Verg\`es forest on $[9]$.

\begin{thm}\label{Josuat-Verges}[Josuat-Verg\`es] For $n\geq 1$, $Q_n(x)$ equals the sum of weights
of all Josuat-Verg\`es forests  on $[n]$.
\end{thm}

By the resemblance between the above combinatorial interpretation of $Q_n(x)$ and
the exponentiation of the Dumont grammar, one sees that the Josuat-Verg\`es forests can be
regarded exactly as an expanded version of the forests for the Dumont grammar.
In our labeling scheme, a leaf is labeled with $u=1+x^2$. This means that
the leaf can be an endpoint (a real leaf) or a vertex having two
empty leaves as its children, and this explains why we should
label an empty leaf by $x$.

We remark that  the exponentiation property
concerning the grammars implies the following
relation on the generating functions of $P_n(x)$ and $Q_n(x)$, as
derived combinatorially with the aid of  cycle alternating permutations by
Josuat-Verg\`es \cite{Josuat-2014},
\begin{equation} \label{g-p-q-exp}
    \sum_{n=0} P_n(x) {t^{n+1} \over (n+1)!} =
    \log \left( \sum_{n=0}^\infty Q_n(x) {t^n \over n!}   \right) =
    \log \left( {1 \over \cos (t) - x \sin(t)} \right),
\end{equation}
which yields the generating function of $P_n(x)$ after differentiation.
For more details about the composition of two grammars, see \cite{Chen-1993}.

\subsection{The $\beta$-expansions}

Ma-Ma-Yeh \cite{Ma-Ma-Yeh-2019} came up with a beautiful idea to deduce the
$\gamma$-positivity of a polynomial by making use of a transformation
of grammars. We will demonstrate that this idea can be adapted
to the derivative polynomials.

 Ma \cite{Ma-2012} obtained  expansions of $P_n(x)$ and $Q_n(x)$
 based on the grammar definitions and a trigonometric identity, which
 we call the $\beta$-expansions. The coefficients
 in the expansions turn out to be the coefficients of the
 peak polynomials. It has been shown that
 \begin{eqnarray}
 P_n(x) & = & \sum_{k=0}^{\lfloor (n-1)/2 \rfloor} M(n,k) x^{n-2k-1}(1+x^2)^{k+1},  \label{Ma-PM}\\[9pt]
 Q_n(x) & = & \sum_{k=0}^{\lfloor n/2 \rfloor} L(n,k)x^{n-2k}(1+x^2)^k. \label{Ma-QT}
  \end{eqnarray}

 As noted by Zhu-Yeh-Lu \cite{Zhu-Yeh-Lu-2019},
 the Springer number $S_n$ equals $L_n(2)$. Similarly, we
 have $P_n(1) = M_n(2)$.
Appealing to the Dumont ansatz for
0-1-2 increasing plane trees,  we see that the $\beta$-expansion of $P_n(x)$
can be viewed as the $\gamma$-expansion of the Eulerian polynomials
$A_n(x,y)$.

\begin{thm}
For $n\geq 1$, we have
\begin{equation}\label{gbeta}
P_n(x) = \sum_{k=1}^{\lfloor (n+1) / 2 \rfloor}
\beta_{n,k} (2x)^{n+1-2k} (1+x^2)^{k},
\end{equation}
where $\beta_{n,k}$ equals the number of 0-1-2 increasing
plane trees on $[n]$ with $k$ leaves.
\end{thm}

Given $u=1+x^2$ and $v=x$, the parameters $(x,y)$ in
(\ref{d-n-a-n}) for the Eulerian polynomials $A_n(x,y)$ are substituted by
\begin{equation}
     (x+i, x-i),
\end{equation}
where $i=\sqrt{-1}$. Hence the following relation is valid.

\begin{thm} For $n\geq 0$, we have
\begin{equation} \label{pnx-a}
    P_n(x) =
    (x-i)^{n+1}  A_n\left( {x+i \over x-i} \right).
\end{equation}
\end{thm}

Again, resorting to the Dumont ansatz, we may associate $P_n(x)$ with
the Andr\'e polynomials.
For $n\geq 1$, we have
\begin{equation}\label{pn-en}
P_n(x) = 2^n E_n\left( {1+x^2 \over 2}, x\right).
\end{equation}

By virtue of
the relation (\ref{E-n-u-v-2}), we may
rewrite (\ref{pn-en}) as
\begin{equation}\label{pn-en-2}
P_n(x) = 2^n x^{n+1} E_n\left( {1+x^2 \over 2x^2}\right).
\end{equation}

For $x=1$, (\ref{pn-en-2}) reduces to
\begin{equation} \label{SE}
    P_n(1)= 2^n E_{n},
\end{equation}
which is due to Knuth-Buckholtz  \cite{Knuth-Buckholtz-1967}.
Notice that the identity (\ref{SE}) also follows from
(\ref{enani}) and (\ref{pnx-a}).
The proof of (\ref{pn-en})
can be regarded as a combinatorial argument.
It is worth mentioning that, as
conjectured by Sun \cite{Sun-2013} and  resolved by Zhu-Yeh-Lu \cite{Zhu-Yeh-Lu-2019},  the sequence of Springer numbers
is log-convex.

The following relation is due to Ma  \cite{Ma-2013}, which specializes to
(\ref{SE}) when $x=0$.
Let $E(x)$ be the generating function of the Euler numbers,
\[ E(x) = \tan (x) + \sec (x).\]
Then
\begin{equation} \label{Ma-identity}
 2^n {{\rm d}^n \over {\rm d}x^n} E(x)  = P_n (E(x)).
\end{equation}

As for the $\beta$-expansion  \eqref{Ma-QT} of $Q_n(x)$,  in light of Theorem \ref{thm-E-D-uv} and Theorem \ref{Josuat-Verges}, one sees that the
coefficient $L(n,k)$ can also be interpreted as the number of forests of planted increasing  binary trees on $[n]$ with $k$ leaves.

Next, we give a grammatical proof of \eqref{Ma-QT}. Recall the grammar for $Q_n(x)$ as given in (\ref{g-q-x}), i.e.,
\begin{equation} \label{g-q-x-2}
G =\{ a \rightarrow a x, \;\; x\rightarrow 1+ x^2 \}.
\end{equation}
Under the substitutions
$v=x$ and $u=1+x^2$,  we get the grammar
\begin{equation} \label{G-avu}
   G = \{ a \rightarrow a v,\;\; v\rightarrow u, \;\; u \rightarrow 2uv\}.
\end{equation}

Let $D$ be the formal derivative of the above grammar $G$. Then we have, for $n \geq 0$,
$$
D^n(a)=a Q_n(x).
$$
The first few  values of $D^n(a)$ are given below,
\begin{eqnarray*}
D(a)&=&a v,\\[6pt]
D^2(a) & = & a(v^2+vu),\\[6pt]
D^3(a)&=& a (v^3+5 v u),    \\[6pt]
D^4(a)&=& a(v^4+18 v^2u+5  u^2),\\[6pt]
D^5(a) & = & a(v^5+58 v^3u+61 vu^2),\\[6pt]
D^6(a) & = & a(v^6+179 v^4u+479 v^2u^2+61 u^3).
\end{eqnarray*}

 \noindent{\it Grammatical Proof of \eqref{Ma-QT}.}
 Setting $z=x^2$,
the grammar
\begin{equation} \label{G-xy-beta}
    G=\{x \rightarrow xy, \;\; y \rightarrow x^2\}
\end{equation}
is transformed into
\begin{equation} \label{G-xyz}
    G=\{ x \rightarrow xy, \;\; y \rightarrow z, \;\; z\rightarrow 2yz \}.
\end{equation}
Now, we may take a comparative look at (\ref{G-avu}) and (\ref{G-xyz}) to
understand what is going on. Recall that for the grammar $G$ in (\ref{G-xy-beta}),
\begin{equation}
    D^n(x) = x \sum_{k=0}^{\lfloor n/2 \rfloor} L(n,k) x^{2k} y^{n-2k}.
\end{equation}
Observe that
the above coefficient $L(n,k)$  of $x^{2k+1} y^{n-2k}$ in $D^n(x)$
 equals the coefficient  of $x  y ^{n-2k} z^{k}$ in $D^n(x)$
 linked with  the grammar in (\ref{G-xyz}), as claimed. \qed

For example,
for the grammar $G$ in (\ref{G-xyz}),  we have
\[
D^5(x)  =  x(y^5 +58  y^{3}z+61 y{z}^{2}).
\]
Renaming $x$ by $a$, $y$ by $x$ and $z$ by $1+x^2$,
the above expression is in accordance with the $\beta$-expansion of $Q_5(x)$:
\[ Q_5(x) = 61 x+180 x^3+120 x^5 = x^5 + 58 x^3 (1+x^2) + 61 x (1+x^2)^2.   \]

To conclude, we remark
that the relation (\ref{Ma-QT})  opens a new avenue for
the Gessel formula for the generating function of $L_n(x,y)$. Yet another
possibility would be to explore the connection between $L_n(x,y)$ and $P_n(x)$
along the line of the exponential formula based on (\ref{g-p-q-exp}).

\vskip 5mm \noindent{\large\bf Acknowledgments.}
We wish to thank S.-M. Ma for stimulating discussions
and the referee for the careful reading and insightful comments.This work was supported by the National
Science Foundation of China.

\end{document}